\documentclass{amsart}

\usepackage{latexsym}
\usepackage{amssymb}

\newtheorem{theorem}{Theorem}[section]

\newtheorem{proposition}[theorem]{Proposition}

\theoremstyle{definition}
\newtheorem{definition}[theorem]{Definition}

\theoremstyle{remark}
\newtheorem{remark}[theorem]{Remark}

\numberwithin{equation}{section}

\newcommand{\C}{\mathbb{C}}
\newcommand{\R}{\mathbb{R}}

\newcommand{\N}{\mathbb{N}}
\newcommand{\E}{\mathbb{E}}

\newcommand{\lm}{\lambda}

\newcommand{\f}{\varphi}

\begin{document}

\title[$C_0$-representations and the Haagerup property]{Notes on $C_0$-representations and the Haagerup property}

\author{Paul Jolissaint}
\address{  Universit\'e de Neuch\^atel,
       Institut de Math\'emathiques,       
       Emile-Argand 11,
       CH-2000 Neuch\^atel, Switzerland}
       
\email{paul.jolissaint@unine.ch}


\date{\today}


\begin{abstract}
For any locally compact group $G$, we show the existence and uniqueness up to quasi-equivalence of a unitary $C_0$-representation $\pi_0$ of $G$ such that the coefficient functions of $C_0$-representations of $G$ are exactly the coefficient functions of $\pi_0$. The present work, strongly influenced by \cite{BG} (which dealt exclusively with discrete groups), leads to new characterizations of
the Haagerup property: $G$ has that property if and only if the representation $\pi_0$ induces a $*$-isomorphism of $C^*(G)$ onto $C^*_{\pi_0}(G)$.
When $G$ is discrete and countable, we also relate the Haagerup property to relative strong mixing properties in the sense of \cite{Jol2} of the group von Neumann algebra $L(G)$ into finite von Neumann algebras.
\end{abstract}

\maketitle

\section{Introduction}

Throughout this article, $G$ denotes a locally compact group. We associate to $G$ a unitary representation $(\pi_0,H_0)$ which has the following properties:
\begin{itemize}
\item it is a $C_0$\textit{-representation}: every coefficient function $s\mapsto \langle\pi_0(s)\xi|\eta\rangle$ associated with $\pi_0$ tends to $0$ as $s\to\infty$;
\item the coefficient functions of $\pi_0$ are exactly the coefficient functions of $C_0$-representations of $G$;
\item the representation $\pi_0$ is the unique $C_0$-representation, up to quasi-equival-ence, which satisfies the above properties.
\end{itemize}
The key idea is to use G. Arsac's notion of $A_\pi$\textit{-spaces} from \cite{Ar}. 

Using the same arguments as in Theorem 3.2 and Corollary 3.4 of \cite{BG}, 
we deduce that:

\par\vspace{3mm}\noindent
\textbf{Proposition A.} 
\emph{Let $G$ be a group as above. Then it has the Haagerup property if and only if the maximal $C^*$-algebra $C^*(G)$ is $*$-isomorphic to the $C^*$-algebra $C^*_{\pi_0}(G)$.}

\par\vspace{3mm}

The preceding proposition deserves a comment which we owe to A. Valette: the Haagerup property of a group $G$ is exactly \textit{property} $C_0$ in the sense of V. Bergelson and J. Rosenblatt in Definition 2.4 of \cite{BR}. Moreover, Theorem 2.5 of the same article states the density of $C_0$-representations in the set of all (classes of) unitary representations on a fixed Hilbert space, and this suffices to prove that there is a $C_0$-representation whose extention to the maximal C$^*$-algebra $C^*(G)$ is faithful.

\par\vspace{3mm}

In the last part of the present notes, we assume that $G$ is discrete and countable. We relate the Haagerup property of $G$ to the embedding of its von Neumann algebra $L(G)$ as a \textit{strongly mixing} subalgebra of some finite von Neumann algebra $M$ in the sense of \cite{Jol2}: this means that, for all $x,y\in M$ such that $\E_{L(G)}(x)=\E_{L(G)}(y)=0$ and for any sequence of unitary operators $(u_n)\subset L(G)$ which converges weakly to $0$, one has
$$
\lim_{n\to\infty}\Vert \E_{L(G)}(xu_ny)\Vert_2=0.
$$
In Section 3, we prove the following result which uses some results from Chapter 2 of \cite{ccjjv}:

\par\vspace{3mm}\noindent
\textbf{Theorem B.}
\emph{Let $G$ be an infinite, countable group. Then it has the Haagerup property if and only if $L(G)$ can be embedded into some finite von Neumann algebra $M$ in such a way that $L(G)$ is strongly mixing in $M$ and that there is a sequence of elements $(x_k)_{k\geq 1}\subset M\ominus L(G)$ such that $\Vert x_k\Vert_2=1$ for every $k$, and
$$
\lim_{k\to\infty}\Vert \lambda(g)x_k-x_k\lambda(g)\Vert_2=0
$$
for every $g\in G$.}

\par\vspace{3mm}\noindent
\textit{Acknowledgements.} We warmly thank A. Valette for his comment about Bergelson and Rosenblatt result mentioned above, and the referee for having detected separability problems in a previous version of the present article and for many valuable comments.

\section{An enveloping $C_0$-representation}

In order to give precise statements of our results, we need to recall some notations and facts on spaces of coefficient functions of unitary representations ($A_\pi$-spaces of G. Arsac) from \cite{Ar} and from P. Eymard's article \cite{Ey}.

The Banach algebra of all continuous functions on $G$ which tend to $0$ at infinity is denoted by $C_0(G)$, and its dense subalgebra formed by all continuous functions with compact support is denoted by $K(G)$. 

Let $(\pi,H)$ be a unitary representation of $G$. If $\xi,\eta\in H$, we denote by
$$
\xi*_\pi \bar\eta(s)=\langle\pi(s)\xi|\eta\rangle\quad (s\in G)
$$
the \textit{coefficient function} associated to $\xi$ and $\eta$. These functions are denoted by $\xi*_\pi\eta$ in \cite{Ar} for instance, but our notation reminds the fact that
$\xi*_\pi \bar\eta$ is linear in $\xi$ and antilinear in $\eta$.

A representation $(\pi,H)$ of $G$ is a \textit{$C_0$-representation} if, for all $\xi,\eta\in H$, the associated coefficient function $\xi*_\pi \bar\eta$ belongs to $C_0(G)$.

The \textit{Fourier-Stieltjes algebra} is the set of all coefficient functions as above. It is denoted by $B(G)$ (\cite{Ey}). 

Recall that $B(G)$ is a Banach algebra with respect to the norm
$$
\Vert \f\Vert_B=\inf\{\Vert\xi\Vert\Vert\eta\Vert : \f=\xi*_\pi\bar{\eta}\}.
$$
It is the dual space of the enveloping $C^*$-algebra $C^*(G)$ under the duality bracket defined on the dense $*$-subalgebra $K(G)$ by
$$
\langle \f,f\rangle=\int_G\f(s)f(s)ds\quad\forall \f\in B(G),\ f\in K(G).
$$

Every unitary representation $(\pi,H)$ of $G$ gives rise to a natural $*$-homomorphism, still denoted by $\pi$, from $C^*(G)$ onto $C^*_\pi(G)$, which extends the map $f\mapsto \pi(f)$ defined on $K(G)$. (Recall that $C^*_\pi(G)$ is the $C^*$-algebra generated by $\{\pi(f):f\in K(G)\}$.)

If $E(G)$ is any subset of $B(G)$, we set
$$
E_1(G)=\{\varphi\in E(G) : \Vert\varphi\Vert_B=1\}
$$
the intersection with the unit sphere of $B(G)$. 

A continuous function $\f:G\rightarrow\C$ is \textit{positive definite} if, for all $s_1,\ldots,s_n\in G$ and all $t_1,\ldots,t_n\in\C$, one has
$$
\sum_{i,j=1}^n\bar{t}_it_j\f(s_i^{-1}s_j)\geq 0.
$$
We denote by $P(G)$ the set of all positive definite functions on $G$. For instance, every coefficient function $\xi*_\pi\bar\xi$ is positive definite, and, conversely, for every $\f\in P(G)$, there exists a unique (up to unitary equivalence) triple $(\pi_\f,H_\f,\xi_\f)$ where $(\pi_\f,H_\f)$ is a unitary representation of $G$ and $\xi_\f$ is a cyclic vector for $\pi_\f$ that satisfies 
$$
\f=\xi_\f*_{\pi_\f}\bar{\xi}_\f.
$$
We recall that $\Vert \f\Vert_B=\f(1)$ for every positive definite function $\f$.

If $\f\in B(G)$, the \textit{adjoint} $\f^*$  of $\f$ is defined by $\f^*(s)=\overline{\f(s^{-1})}$ for every $s\in G$. 
We say that $\f$ is \textit{selfadjoint} if $\f^*=\f$ and we denote by $B_{sa}(G)$ the real Banach algebra of all selfadjoint elements of $B(G)$. Every element $\f\in B_{sa}(G)$ admits a unique decomposition, called \textit{Jordan decomposition}, as 
$$
\f=\f^+-\f^-
$$
where $\f^\pm\in P(G)$ and $\Vert\f\Vert_B=\Vert\f^+\Vert_B+\Vert\f^-\Vert_B$. Thus $B_{sa}(G)=P(G)-P(G)$.

The obvious decomposition of any $\psi\in B(G)$ 
$$
\psi=\frac{1}{2}(\psi+\psi^*)+i\cdot\frac{1}{2i}(\psi-\psi^*)
$$
and the Jordan decomposition imply that 
$$
B(G)=P(G)-P(G)+iP(G)-iP(G).
$$

\par\vspace{3mm}
We also need to recall the definition and a few facts on $A_\pi$-spaces in the sense of G. Arsac \cite{Ar} since they play an important role in the present notes. If $(\pi,H)$ is a unitary representation of $G$, $A_\pi(G)$ is the norm closed subspace of $B(G)$ generated by the coefficient functions $\xi*_\pi\bar\eta$ of $\pi$. Every element $\f\in A_\pi(G)$ can be written as
$$
\varphi= \sum_n \xi_n*_\pi\bar{\eta}_n
$$
where $\xi_n,\eta_n\in H$ for every $n$, $\sum_n\Vert\xi_n\Vert\Vert\eta_n\Vert<\infty$, and where
$$
\Vert\varphi\Vert_B=\inf\{\sum_n\Vert\xi_n\Vert\Vert\eta_n\Vert : \varphi= \sum_n \xi_n*_\pi\bar{\eta}_n\}.
$$ 
The Banach space $A_\pi(G)$ identifies with the predual of the von Neumann algebra $L_\pi(G):=\pi(G)''\subset B(H)$ under the duality bracket
$$
\langle \varphi,\pi(f)\rangle=\int_G \varphi(g)f(g)dg
$$
for every $\f\in A_\pi(G)$ and every $f\in K(G)$.

As is usually the case, $\lambda$ denotes the left regular representation of $G$, and $L(G)=L_\lambda(G)$ is its \textit{associated von Neumann algebra}. In this case, $A(G)=A_\lambda(G)$ is the \textit{Fourier algebra} of $G$ (\cite{Ey}, Chapter 3). 

If $M$ is a von Neumann algebra, its predual is denoted by $M_*$, and if $\f\in M_*$ and $a\in M$, we define $a\f$ and $\f a\in M_*$ by
$$
\langle a\f,x\rangle=\langle \f,xa\rangle\quad\textrm{and}\quad
\langle \f a,x\rangle=\langle \f,ax\rangle\quad\forall x\in M.
$$
Hence, one has $(a_1a_2)\f=a_1(a_2\f)$ and $\f(a_1a_2)=(\f a_1)a_2$ for all $\f\in M_*$ and $a_1,a_2\in M$. 
If $(\pi,H)$ is a unitary representation of $G$, if $\f=\sum_n \xi_n*_\pi\bar{\eta}_n\in A_\pi(G)$, then
$$
\langle\f,x\rangle=\sum_n \langle x\xi_n|\eta_n\rangle\quad\forall x\in L_\pi(G).
$$
If $a\in L_{\pi}(G)$, it is easily checked that
$$
a\f=\sum_n (a\xi_n)*_\pi\bar{\eta}_n\quad \text{and}\quad \f a=\sum_n \xi_n*_\pi\overline{a^*\eta_n}.
$$

Finally, if $(\pi_1,H_1)$ and $(\pi_2,H_2)$ are two unitary representations of $G$, then:
\begin{enumerate}
\item [(1)] we say that they are \textit{quasi-equivalent} if the map $\pi_1(f)\mapsto \pi_2(f)$, from $\pi_1(K(G))$ to $\pi_2(K(G))$, extends to an isomorphism of $L_{\pi_1}(G)$ onto $L_{\pi_2}(G)$;
\item [(2)] we say that they are \textit{disjoint} if no non-zero subrepresentation of $\pi_1$ is equivalent to some subrepresentation of $\pi_2$.
\end{enumerate}

It follows from Propositions 3.1 and 3.12 of \cite{Ar} that:
\begin{enumerate}
\item [(a)] the representations $\pi_1$ and $\pi_2$ are quasi-equivalent if and only if 
$$
A_{\pi_1}(G)=A_{\pi_2}(G);
$$
\item [(b)] the representations $\pi_1$ and $\pi_2$ are disjoint if and only if 
$$
A_{\pi_1}(G)\cap A_{\pi_2}(G)=\{0\}.
$$
\end{enumerate}

\par\vspace{3mm}
Let us now introduce one of the main objects of the present article: let $A_0(G)=B(G)\cap C_0(G)$ be the space of all elements of $B(G)$ that tend to $0$ at infinity. We also put $P_0(G)=P(G)\cap C_0(G)$, and let $A_{0,sa}(G)$ be the real subspace of selfadjoint elements of $A_0(G)$. 

The following result is inspired by \cite{BG}.

\begin{proposition}
The set $A_0(G)$ is a closed two-sided ideal of $B(G)$, it is equal to the set of all coefficient functions of all $C_0$-representations and every $\f\in A_0(G)$ can be expressed as
$$
\f=\f_1-\f_2+i\f_3-i\f_4
$$
with $\f_j\in P_0(G)$ for all $j=1,\ldots,4$.
\end{proposition}
\textit{Proof.} The space $A_0(G)$ is obviously a two-sided ideal of $B(G)$. It is closed because of the following inequality, which holds for every element $\f\in B(G)$:
$$
\Vert\f\Vert_\infty\leq \Vert\f\Vert_B.
$$
Finally, the decomposition of $\f$ as
$$
\varphi=\frac{1}{2}(\varphi+\varphi^*)+i\cdot\frac{1}{2i}(\varphi-\varphi^*)
$$
shows that it suffices to prove that for every selfadjoint element $\f\in A_0(G)$, the positive definite functions $\f^\pm$ of the Jordan decomposition $\f=\f^+-\f^-$ both belong to $C_0(G)$. But it is proved in Lemme 2.12 of \cite{Ey} that $\f^+$ and $\f^-$ are uniform limits on $G$ of linear combinations of right translates $s\mapsto\f(sg)$ of $\f$. As every such translate belongs to $C_0(G)$, this proves the claim.
\hfill $\square$

\par\vspace{3mm}
The reason why we denote the intersection $B(G)\cap C_0(G)$ by $A_0(G)$ instead of $B_0(G)$ for instance is that we will see that it is an $A_\pi$-space for some suitable representation that we introduce now.

We choose some dense directed set $(\f_i)_{i\in I}$ in $P_{0,1}(G)$ and, for every $i\in I$, let $(\pi_i,H_i,\xi_i)$ be the associated cyclic representation. Put first $K_0=\bigoplus_{i\in I} H_i$ and $\sigma_0=\bigoplus_{i\in I}\pi_i$. 
For instance, if $G$ is assumed to be discrete, one can set $\f_1=\delta_1$, so that $\pi_1=\lambda$ is the left regular representation of $G$. Next, set 
$$
H_0=K_0\otimes\ell^2(\N)\quad
\text{and}\quad
\pi_0=\sigma_0\otimes 1_{\ell^2(\N)}.
$$
Notice that both $\sigma_0$ and $\pi_0$ are $C_0$-representations.

\begin{proposition}
Let $G$ be a locally compact, second countable group, and let $(\pi_0,H_0)$ be the above representation. Then:
\begin{enumerate}
\item [(1)] For every $C_0$-representation $\pi$ of $G$, one has $A_\pi(G)\subset A_0(G)$.
\item [(2)] One has $A_0(G)=A_{\pi_0}(G)$, and every coefficient function of any $C_0$-repre-sentation is a coefficient function associated to $\pi_0$.
\item [(3)] The unitary representation $\pi_0$ is the unique $C_0$-representation such that $A_0(G)=A_{\pi_0}(G)$, up to quasi-equivalence.
\end{enumerate}
\end{proposition}
\textit{Proof.} (1) Observe that every coefficient function $\f$ of the $C_0$-representation $\pi$ is a linear combination of four elements in $P_{0,1}(G)$, by the same argument as in the proof of Proposition 2.1. As $A_0(G)$ is closed, this proves the first assertion. In particular, $A_{\sigma_0}(G)$ and $A_{\pi_0}(G)$ are contained in $A_0(G)$.\\
(2) First, if $\f\in P_{0,1}(G)$, then it is a norm limit of a subsequence $(\psi_{k})_{k\geq 1}$ of $(\f_i)$. This shows that $\f\in A_{\sigma_0}(G)$, and Proposition 2.1 proves that $A_0(G)\subset A_{\sigma_0}(G)\subset A_{\pi_0}(G)$.
Next, let $\f\in A_0(G)$. Let us prove that it is a coefficient function of $\pi_0$. As $A_{\sigma_0}(G)=A_0(G)$, there exist sequences of vectors $(\xi_n)_{n\geq 1},(\eta_n)_{n\geq 1}\subset K_0$ such that 
$$
\sum_n\Vert\xi_n\Vert\Vert\eta_n\Vert<\infty
$$
and 
$$
\f=\sum_n\xi_n*_{\sigma_0}\bar{\eta}_n.
$$
Replacing $\xi_n$ by $\sqrt{\frac{\Vert\eta_n\Vert}{\Vert\xi_n\Vert}}\xi_n$ and $\eta_n$ by
$\sqrt{\frac{\Vert\xi_n\Vert}{\Vert\eta_n\Vert}}\eta_n$, we assume that 
$$
\sum_n\Vert\xi_n\Vert^2=\sum_n\Vert\eta_n\Vert^2=\sum_n\Vert\xi_n\Vert\Vert\eta_n\Vert<\infty.
$$
Put $\xi=\bigoplus_n \xi_n,\eta=\bigoplus_n\eta_n\in H_0$. Then $\f=\xi*_{\pi_0}\bar\eta$.\\
(3) follows immediately from (1) and (2).
\hfill $\square$

\begin{definition}
The representation $(\pi_0,H_0)$ is called the \textbf{enveloping $C_0$-repre-sentation} of $G$.
\end{definition}

\begin{remark}
(1) As is well known, the left regular representation of $G$ is a $C_0$-representation. Hence the Fourier algebra $A(G)$ is contained in $A_0(G)$. In fact, one can have equality $A(G)=A_0(G)$ as well as strict inclusion 
$A(G)\subsetneq A_0(G)$. Indeed, on the one hand, I. Khalil proved in \cite{Kha} that if $G$ is the $ax+b$-group over $\R$, then $A(G)=A_0(G)$, and, on the other hand, A. Fig\`a-Talamanca \cite{FG} proved that if $G$ is unimodular and if its von Neumann algebra $L(G)$ is not atomic (e.g. it is the case whenever $G$ is infinite and discrete), then $A(G)\subsetneq A_0(G)$.\\
(2) We are grateful to the referee for the following observation: the proofs of Propositions 2.1 and 2.2 show that they hold with $A_0(G)$ replaced by any norm-closed, $G$-invariant subspace of $B(G)$.
\end{remark}

The next proposition is strongly inspired by, and is a slight generalization of Theorem 3.2 of \cite{BG}. It will be used to give characterizations of the Haagerup property in terms of the enveloping $C_0$-representation.

\begin{proposition}
Let $G$ be locally compact group and let $(\pi,H)$ be a unitary representation of $G$, and let us assume that the space $A_\pi(G)$ is an ideal of $B(G)$. Then $\pi:C^*(G)\rightarrow C^*_\pi(G)$ is a $*$-isomorphism if and only if there is a sequence $(\f_n)_{n\geq 1}\subset A_\pi(G)\cap  P_1(G)$ such that $\f_n\to 1$ uniformly on compact subsets of $G$.
\end{proposition}
\textit{Proof.} Assume first that $\pi$ is a $*$-isomorphism. We can suppose that $C^*_\pi(G)$ contains no non-zero compact operator. Let $\chi$ be the state on $C^*_\pi(G)$ which comes from the trivial character $f\mapsto \int_G f(s)ds$ on $K(G)\subset C^*(G)$. By Glimm's Lemma, there exists an orthonormal sequence $(\xi_n)_{n\geq 1}\subset H$ such that 
$$
\chi(x)=\lim_{n\to\infty}\langle x\xi_n|\xi_n\rangle
$$
for every $x\in C^*_\pi(G)$. Put $\f_n=\xi_n*_{\pi}\bar{\xi}_n\in A_\pi(G)\cap P_1(G)$ for every $n$. Then one has for every $f\in K(G)$:
$$
\lim_{n\to\infty}\int_G\f_n(t)f(t)dt=\lim_{n\to\infty}\langle\pi(f)\xi_n|\xi_n\rangle=\int_G f(t)dt.
$$
Theorem 13.5.2 of \cite{DiC} implies that $\f_n\to 1$ uniformly on compact subsets of $G$.
\\
Conversely, if there exists a sequence $(\f_n)_{n\geq 1}\subset A_\pi(G)\cap P_1(G)$ such that $\f_n\to 1$ uniformly on compact subsets of $G$, let $x\in\mathrm{ker}(\pi)$. We have to prove that $\langle\f,x^*x\rangle_{B,C^*}=0$ for every state $\f$ on $C^*(G)$. Observe first that, for every $\psi\in A_\pi(G)$ and every $y\in C^*(G)$, one has 
$$
\langle \psi,y\rangle_{B,C^*}=\langle \psi,\pi(y)\rangle_{A_\pi,C^*_\pi}.
$$
Indeed, if we write $\psi=\sum_k\xi_k*_\pi\bar{\eta}_k$, and if $f\in K(G)$, we have
$$
\langle \psi,f\rangle_{B,C^*}=\int_G\psi(s)f(s)ds=
\sum_k \int_G\langle\pi(s)\xi_k|\eta_k\rangle f(s) ds=\langle \psi,\pi(f)\rangle_{A_\pi,C^*_\pi}
$$
and the formula holds by density of $K(G)$ in $C^*(G)$.

Let us fix such a state $\f\in P_1(G)$ and set $\psi_n=\f\f_n\in A_\pi(G)\cap P_1(G)$ for every $n$. As $\psi_n$ is a state on $L_\pi(G)$, its restriction to $C^*_\pi(G)$ is still a state, and 
$\langle \psi_n,x^*x\rangle=\langle \psi_n,\pi(x^*x)\rangle=0$ for every $n$. As $\psi_n\to\f$ in the weak$^*$ topology of $B(G)=C^*(G)^*$, one has $\langle \f,x^*x\rangle=0$.
\hfill $\square$

\section{The Haagerup property}

As in the first section, $G$ denotes a locally compact group and $(\pi_0,H_0)$ denotes its enveloping $C_0$-representation.

Following M. Bekka \cite{Bekka}, we say that $(\pi,H)$ is an \textit{amenable representation} if $\pi\otimes\bar\pi$ weakly contains the trivial representation. Equivalently, this means that there exists a net of unit vectors $(\xi_i)\subset H\otimes \bar H$ such that 
$$
\langle\pi\otimes\bar\pi(s)\xi_i|\xi_i\rangle\to 1
$$
uniformly on compact subsets of $G$; notice that $\pi\otimes\bar\pi$ is unitarily equivalent to the representation $(T,g)\mapsto \pi(g)T\pi(g^{-1})$ acting on the space $HS(H)$ of all Hilbert-Schmidt operators.

If $G$ is moreover second countable, we say that it has the \textit{Haagerup property} if there exists a sequence $(\f_n)_{n\geq 1}\subset P_{0,1}(G)$ which tends to $1$ uniformly on compact sets. Note that it is equivalent to say that $G$ admits an amenable, $C_0$-representation. See \cite{ccjjv} for more information on the Haagerup property.

The next result generalizes partly, and is inspired by Corollary 3.4 of \cite{BG}. 

\begin{proposition}
Let $G$ and $(\pi_0,H_0)$ be as above. Then the following conditions are equivalent:
\begin{enumerate}
\item [(1)] $G$ has the Haagerup property;
\item [(2)] $C^*(G)=C^*_{\pi_0}(G)$, i.e. the $*$-homomorphism $\pi_0:C^*(G)\rightarrow C^*_{\pi_0}(G)$ is an isomorphism;
\item [(3)] the representation $\pi_0$ weakly contains the trivial representation;
\item [(4)] the representation $\pi_0$ is amenable.
\end{enumerate}
\end{proposition}
\textit{Proof.} $(1)\ \Rightarrow\ (2)$. There exists a sequence $(\f_n)_{n\geq 1}\subset P_{0,1}(G)$ which converges to $1$ uniformly on compact sets. The assertion follows readily from Proposition 2.5.\\
$(2)\ \Rightarrow\ (3)$. It follows also from Proposition 2.5.\\
$(3)\ \Rightarrow\ (4)$ and $(4)\ \Rightarrow\ (1)$ are obvious.
\hfill $\square$

\begin{remark}
As $A(G)\subset A_{\pi_0}(G)$, there exists a $*$-homomorphism $\Phi$ from $L_{\pi_0}(G)$ onto $L(G)$ such that $\Phi(\pi_0(f))=\lambda(f)$ for every $f\in K(G)$. Thus, let $z_A\in L_{\pi_0}(G)$ be the central projection such that  $L_{\pi_0}(G)z_A$ is $*$-isomorphic to $L(G)$. This allows us to consider the following two subrepresentations of $\pi_0$: set $\pi_{00}(s)=\pi_0(s)(1-z_A)$ and $\lambda_0(s)=\pi_0(s)z_A$ for all $s\in G$. Then $\lambda_0$ is quasi-equivalent to $\lambda$, and since $\pi_{00}$ is disjoint from $\lambda$, we have $A_{\pi_{00}}(G)\cap A(G)=\{0\}$. It would be interesting to get more information on $\pi_{00}$, in particular when $G$ has the Haagerup property.
\end{remark}

\par\vspace{3mm}
From now on, we assume that $G$ is an infinite, discrete, countable group. Following \cite{BG}, for any (not necessarily closed) ideal $D\subset \ell^\infty(G)$, we say that a unitary representation $(\pi,H)$ of $G$ is a \textit{$D$-representation} if $H$ contains a dense subspace $K$ such that the coefficient function $\xi*_\pi\bar\eta\in D$ for all $\xi,\eta\in K$. We associate to $D$ the following $C^*$-algebra $C^*_D(G)$: it is the completion of $K(G)$ with respect to the $C^*$-norm
$$
\Vert f\Vert_D:=\sup\{\Vert\pi(f)\Vert : \pi\ \text{is\ a}\ D-\text{representation}\}.
$$
When $D=C_0(G)$, one gets $C^*_D(G)=C^*_{\pi_0}(G)$. This makes the link between Proposition 3.1 above and the main results of N. Brown and E. Guentner in \cite{BG}.

\par\vspace{3mm}
We end the present notes with a relationship between the Haagerup property for discrete groups and strongly mixing von Neumann subalgebras in the sense of \cite{Jol2}, Definition 1.1. We need  to recall some definitions and facts from \cite{Jol2} first and from Chapter 2 of \cite{ccjjv} next.

Let $1\in B\subset M$ be finite von Neumann algebras (with separable preduals) endowed with a normal, finite, faithful, normalized trace $\tau$. We denote by $\E_B$ the $\tau$-preserving conditional expectation from $M$ onto $B$, and by $M\ominus B=\{x\in M:\E_B(x)=0\}$. We assume that $B$ is diffuse.

\begin{definition}
Let $B\subset M$ be a pair as above. 
We say that $B$ is \textbf{strongly mixing in} $M$ if 
$$
\lim_{n\to\infty}\Vert \E_B(xu_ny)\Vert_2=0
$$
for all $x,y\in M\ominus B$ and all sequences $(u_n)\subset U(B)$ which converge to $0$ in the weak operator topology.
\end{definition}

 This definition is motivated by the following situation: if a countable group $G$ acts in a trace-preserving way on some finite von Neumann algebra $(Q,\tau)$ and if we put $B:=L(G)\subset M:=Q\rtimes G$, then $B$ is strongly mixing in $M$ if and only if the action of $G$ on $Q$ is strongly mixing in the usual sense: for all $a,b\in Q$, one has $
\lim_{g\to\infty}\tau(a\sigma_g(b))=\tau(a)\tau(b).$

Let now $G$ be a countable group with the Haagerup property. By Theorems 2.1.5, 2.2.2 and 2.3.4 of \cite{ccjjv}, there exists a trace preserving and strongly mixing action of $G$ on some finite von Neumann algebra $(Q,\tau)$ which contains non trivial asymptotically invariant sequences and F\o lner sequences in the sense below. For instance, if $G$ has the Haagerup property, there exists an action $\alpha$ of $G$ on the hyperfinite type II$_1$-factor $R$ such that:
\begin{itemize}
\item $\alpha$ is strongly mixing;
\item the fixed point algebra $(R_\omega)^\alpha$, that is, the set of all (classes of) central sequences $x=[(x_n)]\in R_\omega$ such that $\alpha_g^\omega(x)=x$ for all $g\in G$, is of type II$_1$. 
\end{itemize}

\begin{definition}
Let $1\in B\subset M$ be a pair of finite von Neumann algebras as above, and let $(e_k)_{k\geq 1}\subset M$ be a sequence of projections in $M$.
\begin{enumerate}
\item [(1)] We say that $(e_k)_{k\geq 1}$ is a \textbf{non trivial asymptotically invariant sequence} for $B$ if $\E_B(e_k)=\tau(e_k)$ for every $k$, if
$$
\lim_{k\to\infty} \Vert be_k-e_kb\Vert_2=0
$$
for every $b\in B$ and if
$$
\inf_k\tau(e_k)(1-\tau(e_k))>0.
$$
\item [(2)] We say that $(e_k)_{k\geq 1}$ is a \textbf{F\o lner sequence} for $B$ if $\E_B(e_k)=\tau(e_k)$ for every $k$, if $\lim_k\Vert e_k\Vert_2=0$ and if
$$
\lim_{k\to\infty} \frac{\Vert be_k-e_kb\Vert_2}{\Vert e_k\Vert_2}=0
$$
for every $b\in B$.
\end{enumerate}
\end{definition}

In general, the existence of a non trivial asymptotically invariant sequence for $B$ implies the existence of a F\o lner sequence for $B$, but the converse does not hold. See \cite{ccjjv}, p. 19, for more details.

\par\vspace{3mm}

Combining these types of properties, we get:

\begin{theorem}
Let $G$ be an infinite, countable group. Then it has the Haagerup property if and only if it satisfies one of the following equivalent conditions:
\begin{enumerate}
\item [(1)] $\textrm{(resp.\ }$ $(1')$) There exists a finite von Neumann algebra $M$ containing $L(G)$ such that $L(G)$ is strongly mixing in $M$ and $M$ contains a F\o lner sequence for $L(G)$ (resp. a non trivial asymptotically invariant sequence for $L(G)$).
\item [(2)] There exists a finite von Neumann algebra $M$ containing $L(G)$ such that $L(G)$ is strongly mixing in $M$ and there is a sequence of elements $(x_k)_{k\geq 1}\subset M\ominus B$ such that $\Vert x_k\Vert_2=1$ for every $k$, and
$$
\lim_{k\to\infty}\Vert \lambda(g)x_k-x_k\lambda(g)\Vert_2=0
$$
for every $g\in G$. 
\end{enumerate}
\end{theorem}
\textit{Proof.} If $G$ has the Haagerup property, then each condition (1), (1') and (2) holds, by Theorem 2.3.4 of \cite{ccjjv}, and there are plenty of non trivial asymptotically invariant or F\o lner sequences in the hyperfinite type II$_1$-factor $R$. Thus, assume that condition (1) holds and that $B:=L(G)$ embeds into some finite von Neumann algebra $M$ such that $B:=L(G)$ is strongly mixing in $M$ and that $M$ contains a F\o lner sequence for $B$. We have to show the existence of a sequence $(\f_k)_{k\geq 1}\subset P_{0,1}(G)$ which tends to $1$ pointwise. 

Recall first that to any completely positive map $\Phi:M\rightarrow M$, one associates a function $\f$ on $G$ by
$$
\f(g)=\tau(\Phi(\lm(g))\lm(g^{-1}))\quad (g\in G),
$$
and that $\f$ is positive definite. In particular, for every $x\in M\ominus B$, the function $\f_x:G\rightarrow \C$ defined by
$$
\f_x(g)=\tau(\E_B(x^*\lm(g)x)\lm(g^{-1}))=\tau(x^*\lm(g)x\lm(g^{-1}))\quad (g\in G)
$$
is positive definite. Moreover, since $B$ is strongly mixing in $M$ and since $\lm(G)$ is an orthonormal set, one has
$$
|\f_x(g)|\leq \Vert\E_B(x^*\lm(g)x)\Vert_2\to 0
$$
as $g\to\infty$, which shows that $\f_x\in P_0(G)$ for every $x\in M\otimes B$. \\
Next, let $(e_k)_{k\geq 1}\subset M$ be a F\o lner sequence for $B$ and choose $c>0$ and an integer $k_0>0$ such that
$$
1-\tau(e_k)\geq c
$$
holds for every $k\geq k_0$.
Define then
$$
x_k=\frac{e_k-\tau(e_k)}{\sqrt{\tau(e_k)(1-\tau(e_k))}}(=x_k^*)\quad (k\geq 1)
$$
and put $\f_k=\f_{x_k}$ for every $k$. One has, for every integer $k\geq k_0$ and every $g\in G$:
\begin{eqnarray*}
\f_k(g) &=&
\tau(x_k\lm(g)x_k\lm(g^{-1}))\\
&=&
\frac{1}{\tau(e_k)(1-\tau(e_k))}\cdot \tau((e_k-\tau(e_k))\lm(g)(e_k-\tau(e_k))\lm(g^{-1}))\\
&=&
\frac{1}{\tau(e_k)(1-\tau(e_k))}\cdot\tau(e_k\lm(g)e_k\lm(g^{-1})-\tau(e_k)^2)\\
&=&
\frac{\tau(e_k(\lm(g)e_k\lm(g^{-1})-e_k))}{\tau(e_k)(1-\tau(e_k))}+1.
\end{eqnarray*}
Hence, by Cauchy-Schwarz Inequality,
\begin{eqnarray*}
|\f_k(g)-1|
&\leq &
\frac{1}{c}\cdot\frac{\Vert e_k\Vert_2\Vert\lm(g)e_k\lm(g^{-1})-e_k\Vert_2}{\Vert e_k\Vert_2^2}\\
& =&
\frac{1}{c}\cdot\frac{\Vert\lm(g)e_k-e_k\lm(g)\Vert_2}{\Vert e_k\Vert_2}
\to 0
\end{eqnarray*}
as $k\to\infty$ for every $g\in G$. 
A similar argument works if $(e_k)$ is a non trivial asympotically invariant sequence.

Finally, assume that $G$ satisfies condition (2), and let $(x_k)\subset M\ominus B$ be as above. Define
$\f_k(g)=\tau(x_k^*\lambda(g)x_k\lambda(g^{-1}))$ exactly as above. Then by the same arguments, $\f_k\in P_{0,1}(G)$ for every $k$, and, for fixed $g\in G$, one has:
\begin{eqnarray*}
|\f_k(g)-1| &= &
|\tau(x_k^*\lm(g)x_k\lm(g^{-1}))-\tau(x_k^*x_k)|\\
&=&
|\langle \lm(g)x_k\lm(g)-x_k|x_k\rangle|\\
&\leq &
\Vert\lm(g)x_k\lm(g^{-1})-x_k\Vert_2\Vert x_k\Vert_2\\
&=&
\Vert\lm(g)x_k\lm(g^{-1})-x_k\Vert_2\to 0
\end{eqnarray*}
as $k\to \infty$.
\hfill $\square$

\begin{remark}
Assume that $G$ has the Haagerup property. One can ask whether there exists a group $\Gamma$ containing $G$ and such that the pair of finite von Neumann algebras $L(G)\subset L(\Gamma)$ satisfies condition (2) in Theorem 3.5. Unfortunately, it is only the case when $G$ is amenable, and this has no real interest. Indeed, assume for simplicity that $G$ is torsion free, that it embeds into some group $\Gamma$ and that the pair $L(G)\subset L(\Gamma)$ satisfies condition (2) above. Then, on the one hand, by Lemma 2.2 and Proposition 2.3 of \cite{Jol2}, the pair of groups $G\subset \Gamma$ satisfies \textit{condition (ST)}, which means that, for every $\gamma\in \Gamma\setminus G$, the subgroup $\gamma G\gamma^{-1}\cap G$ is finite, hence trivial. In other words, $G$ is \textit{malnormal} in $\Gamma$. On the other hand, by classical arguments, the existence of a sequence $(x_k)\subset L(\Gamma)\ominus L(G)$ as above implies that the action $G\curvearrowright X:=\Gamma\setminus G$ defined by $(g,x)\mapsto gxg^{-1}$ has an invariant mean.
This means that the associated representation $\lambda_X$ weakly contains the trivial representation. But the first condition implies that this action is free, hence that $\lambda_X$ is equivalent to a multiple of the regular representation. This forces $G$ to be amenable. 
\end{remark}

\par\vspace{3mm}

\bibliographystyle{plain}
\bibliography{ref}

\end{document}